\documentclass{amsart}

\usepackage{amsthm, amsmath,amssymb,stmaryrd}
\usepackage{maple2e}
\usepackage[all]{xy}

\usepackage[pdfauthor   = {Arne Lorenz},
            pdftitle    = {On Local Integrability Of Jet Groupoids},
            pdfsubject  = {58H05 \and 58A20 \and 58A32},
            pdfkeywords = {},
            bookmarks=true,
            bookmarksopen=true,
            colorlinks=true,
            pagebackref=true,
            hyperindex=true,
            linkcolor=blue,
            citecolor=blue,
            urlcolor=blue,
            hypertex=true
            ]{hyperref}

\numberwithin{equation}{section}

\newcommand{\N}{{\mathbb{N}}}
\newcommand{\Z}{{\mathbb{Z}}}
\newcommand{\RR}{{\mathbb{R}}}

\newcommand{\mcE}{\mathcal{E}}
\newcommand{\mcF}{\mathcal{F}}
\newcommand{\mcR}{\mathcal{R}}

\DeclareMathOperator{\img}{im}
\DeclareMathOperator{\Stab}{Stab}

\DeclareMathOperator{\pr}{pr}

\DeclareMathOperator{\GL}{GL}

\newcommand\fibreProd[3]{\mathinner{ \underset{#2}{{}^{#1} \!
       \Yup  \!  {}^{#3}}}}

\newtheorem{theorem}{Theorem}[section]

\theoremstyle{definition}
\newtheorem{definition}[theorem]{Definition}

\newtheorem{proposition}[theorem]{Proposition}

\begin{document}

\begin{abstract}
  A Jet groupoid $\mcR_q$ over a manifold $X$ is a special Lie groupoid
  consisting of $q$-jets of local diffeomorphisms $X \to X$. As a subbundle of
  $J_q(X \times X)$, a jet groupoid can be considered as a nonlinear system of
  partial differential equations (PDE). This leads to the concept of formal
  integrability. On the other hand, each jet groupoid is the symmetry groupoid
  of a geometric object, modelled as a section $\omega$ of a natural bundle
  $\mcF$. Using the jet groupoids, we give a local characterisation of formal
  integrability for transitive jet groupoids in terms of their corresponding
  geometric objects.
\end{abstract}

\title{On Local Integrability Conditions Of Jet Groupoids} 

\author{Arne Lorenz}
\address{Lehrstuhl B f\"ur Mathematik, RWTH Aachen University, 52062 Germany} 
\email{\href{mailto:Arne Lorenz <arne@momo.math.rwth-aachen.de>}{arne@momo.math.rwth-aachen.de}}

\date{\today}


\thanks{Thanks to M. Barakat for discussions. The author was supported by DFG
    Grant Graduier\-ten\-kol\-leg 775.}
\maketitle




\section{Introduction}
\label{sec_Introduction}

In two articles \cite{Lie1891}, Lie introduced pseudogroups and formulated
their defining equations as differential invarants of the corresponding
pseudogroup. Vessiot continued in \cite{Vessiot03} with the calculation of
necessary conditions for integrability. Later, Pommaret \cite{Pommaret78}
applied Spencer's approach \cite{Spencer69} to pseudogroups and showed that they
are all given as symmetry transformations of geometric objects $\omega$ on
natural bundles $\mcF$.

In this paper, we use the language of jet groupoids which provides an efficient
language for this theory, including conceptional proofs. A Jet groupoid
consists of the algebraic solutions of the pseudogroup equations and is thus also
defined by a section of a natural bundle. We express the prolongation and projection
of jet groupoids by sections of new natural bundles. All considerations are local. As the main
result, we obtain conditions on the sections $\omega$ of $\mcF$ that classify
all formally integrable jet groupoids that can be defined by sections of $\mcF$. 

The last section gives explicit calculations for the well-known example of a
Riemannian metric on a two-dimensional manifold using computer algebra.

\section{Preliminaries}
\label{sec_Preliminaries}

Before turning to jet groupoids, we shortly introduce the underlying Lie
grou\-poids and their action on fibre bundles. This is done to fix the
notation, recent introductory books on Lie groupoids are \cite{Mackenzie05} and
\cite{MoerMrcun03}, which also provide many further references.

\begin{definition}
  \label{def_Groupoid}
  \index{Groupoid}
  A groupoid $G$ is a small category with invertible morphisms. 
\end{definition}

The set of objects, denoted by \index{$G^{(0)}$} $G^{(0)}$, is called the \emph{base}. The
set of morphisms is denoted by \index{$G^{(1)}$} $G^{(1)}$ and has the
projections \emph{source} and \emph{target} to the base:
\[
s : G^{(1)} \to G^{(0)}:  g \mapsto x \qquad t : G^{(1)} \to G^{(0)}:  g \mapsto
y \qquad \mbox{for } g : x \to y \in G^{(1)}.
\]
Composition of morphisms induces a partial multiplication that is defined
whenever source and target match:
\[
\mu : 
G^{(1)} \fibreProd{s \!}{G^{(0)}}{\! t}  G^{(1)} \to G^{(1)}: (g,h) \mapsto gh 
\]
with $G^{(1)} \, {}^{s \! \!} \! \Yup^{\! t}_{G^{(0)}}  G^{(1)} = \{ (g,h) \in G^{(1)} \times G^{(1)} |
s(g) = t(h) \}$. Via $\iota: G^{(0)} \hookrightarrow G^{(1)}: x \mapsto 1_x$,
the base is embedded in $G^{(1)}$ as the identity morphisms. A groupoid is called
\emph{transitive} if $G(x,y) = \{g \in G^{(1)} | s(g)=x, t(g)=y \} \neq
\emptyset$ for all $x,y \in G^{(0)}$.

\begin{definition}
  \label{def_LieGroupoid}
  \index{Groupoid!Lie}
  \index{Lie Groupoid}
  A groupoid $G$ is called a \emph{Lie groupoid} if $G^{(0)}$ and $G^{(1)}$ are
  smooth\footnote{In this paper, smooth means $C^\infty$.}
  manifolds where $s$, $t$, $\mu$, $\iota$ and the inversion are smooth, $s$,
  $t$ being surjective submersions.
\end{definition}

In the case of Lie groupoids, $G^{(1)} \,{}^{s \! \!} \! \Yup^{\! t}_{G^{(0)}}
G^{(1)}$ turns into the fibre product and the isotropy groups $G(x,x)$
are Lie groups \cite[Thm 5.4]{MoerMrcun03}. Important examples of Lie
groupoids are gauge groupoids $PP^{-1} := P \times_H P$ for Lie groups $H$ and
principal $H$-bundles $P$.

\begin{definition}
  \label{def_GroupoidAction}
  \index{Groupoid!action}
  Let $G$ be a Lie groupoid and $\pi: \mathcal{F} \to G^{(0)}$ be a fibre
  bundle. A right \emph{groupoid action} of $G$ on $\mcF$ is a smooth map $\mcF
  \, {}^{\pi \! \!} \! \Yup^{\! t}_{G^{(0)}}  G^{(1)} \to \mcF$ with $f(ab) = (fa)b$ and
  $f 1_x = f$ whenever $f \in \mcF_x = \pi^{-1}(x)$ and $a,b \in G^{(1)}$ can be composed.
\end{definition}

\section{Jet Groupoids and Natural Bundles}
\label{sec_Jet_Groupoids}

We will now define jet groupoids as Lie groupoids over a fixed base manifold
$G^{(0)} = X$ of dimension $n$. The $q$-th jet bundle $J_q(X \times X)$ over the
trivial bundle $X \times X$ provides source $s = \pr_1$ and target $t = \pr_2$
as projections on the first and second copy of $X$. Restricted to the open
subset $\Pi_q \subset J_q(X \times X)$ of invertible jets, the chain rule
induces a partial multiplication on $\Pi_q$, which is respected by the natural
projections $\pi^{q+r}_q : J_{q+r}(X \times X) \to J_q(X \times X)$.

\begin{definition}
  \label{def_JetGroupoid}
  $\Pi_q$ is called the \emph{full jet groupoid} of order $q$ and a \emph{jet
    groupoid} $\mcR_q$ is a subbundle of $\Pi_q$, closed with respect to all
  groupoid operations. 
\end{definition}

For the treatment of natural bundles and the projection of jet groupoids, it is
helpful to consider the isotropy groups $\Pi_q(x,x)$. They are all isomorphic to
$\GL_q = \GL_q(\RR^n)$, the Lie group of $q$-jets of diffeomorphisms
$\RR^n \to \RR^n$ leaving the origin fixed. The structure of $\GL_q$ was studied
by Terng \cite{Terng78}. By the construction of $\GL_{q+1}$ there is an exact sequence
\begin{equation}
  \label{eq_GL_exact_seq}
  \xymatrix{
    1 \ar[r] & K_{q+1} \ar[r] & \GL_{q+1} (\RR^n) \ar[r]^{\pi^{q+1}_q} & \GL_q
    (\RR^n) \ar[r] & 1
  }
\end{equation}
defining the normal subgroup $K_{q+1} := \ker(\pi^{q+1}_q) \unlhd
\GL_{q+1}$. The projection can be identified with $\pi^{q+1}_q : \Pi_{q+1}(x,x)
\to \Pi_q(x,x)$ for each $x \in X$. $\GL_{q+1}$ is called first principal
prolongation of $\GL_q$ in \cite{KMS93} (following \cite{Ehresmann55}).

All jets with fixed target $y_0 \in X$ define a
principal $\GL_q$-bundle $P_q := \Pi_q(-,y_0)$. When changing from groupoids to
the bundle point of view, left and right $\GL_q \cong \Pi_q(y_0,y_0)$-actions
must be swapped to obtain the equations in \cite{Lie1891}, \cite{Vessiot03} and
(right) principal bundles. We recover $\Pi_q$ as the gauge groupoid $P_q
P_q^{-1}$ (via $(g,h) \mapsto g h^{-1}$). The sequence (\ref{eq_GL_exact_seq})
implies $P_{q+1}/ K_{q+1} \cong P_q$. Writing $K_{q+1}$ for all kernels
$\ker(\pi^{q+1}_q) \unlhd \Pi_{q+1}(x,x)$, we obtain a kind of commutation law
$K_{q+1} f_{q+1} = f_{q+1} K_{q+1}$ as sets for all $f_{q+1} \in \Pi_{q+1}$.

\begin{definition}
  \label{def_NaturalBundle}
  A fibre bundle $\mathcal{F} \stackrel{\pi}{\rightarrow} X$ is called
  \emph{natural bundle} if there exists a $q \in \N$ and a groupoid action of
  $\Pi_q$ on $\mathcal{F}$. A section $\omega$ of $\mcF$ is called \emph{geometric
  object}.
\end{definition}

$P_q$ is a natural bundle by right $\Pi_q$-multiplication. In fact, all natural
bundles $\mcF$ with typical fibre $F := \mcF_{y_0}$ are associated to $P_q$ as
$\mcF \cong P_q \times_{\GL_q} F$. This is done by splitting $u \in \mcF$ into
$u = u_{y_0} f_q$ with $u_{y_0} \in F$ and $f_q \in P_q$, unique up to elements
of $\GL_q \cong \Pi_q(y_0,y_0)$.

If not stated otherwise, we assume the natural bundles $\mcF$ to have fibres $F$
that are homogeneous $\GL_q$-spaces. 

\begin{proposition}
  \label{prop_Groupoid_Bundle}
  Each section $\omega$ of a natural bundle $\mcF$ defines a jet groupoid
  $\mcR_q(\omega) = \Stab_\mcF^q(\omega)$. Conversely, each transitive jet
  groupoid $\mcR_q$ defines a natural bundle $\mcF$ with
  section $\omega_0$, such that $\mcR_q$ is the full
  symmetry groupoid $\Stab_\mcF^q(\omega_0)$ of $\omega_0$.
\end{proposition}

\begin{proof}
  Define the symmetry groupoid $\Stab_\mcF^q(\omega)$ via the $\Pi_q$-action
  on $\mcF$
  \[
  \Phi_\omega : \Pi_q \to \mcF : f_q \mapsto \omega(t(f_q)) f_q
  \]
  as the kernel $\ker_\omega(\Phi_\omega) =
  \{f_q \in \Pi_q | \Phi_\omega(f_q) = \omega(s(f_q)) \}$ or by the exact
  sequence
  \begin{equation}
    \label{eq_R_q_sequence}
    \xymatrix{
      0 \ar[r] & 
      \Stab_\mcF^q(\omega)
      \ar[r] & \Pi_q \ar@<0.6ex>[r]^{\Phi_\omega }
      \ar@<-0.4ex>[r]_{\omega \circ s}  & \mcF
    }
  \end{equation}
  of bundles over $X = s(\Pi_q)$. The $\Pi_q$-action implies that
  $\Stab_\mcF^q(\omega)$ is a groupoid since it is closed under $\mu$, $\iota$
  and inversion. As $F$ is homogeneous, $\Phi_\omega$ is surjective and of
  constant rank. By the implicit function theorem, $\Stab_\mcF^q(\omega)$ is a
  Lie groupoid. Each $f_q \in \Pi_q$ can be modified by $g_q \in \GL_q$ such
  that $\omega(y) f_q g_q = \omega(x)$, so $\Stab_\mcF^q(\omega)$ is transitive.
  
  The transitivity of $\mcR_q$ implies that all isotropy groups are isomorphic
  to some $G_q \leq \GL_q$, so choose $y_0 \in X$ and set $\mcF := \mcR_q(y_0,
  y_0) \backslash \Pi_q(-,y_0) \cong P_q \times_{\GL_q} \GL_q/G_q$ with the section
  \[
  \omega_0 : X \rightarrow \mcF : x \mapsto \mcR_q(y_0,y_0) 
  \mcR_q(x,y_0) = \mcR_q(x,y_0).
  \]
  The condition for $r_q \in \Stab_\mcF^q(\omega_0)$ is $\omega_0(y) r_q =
  \omega_0(y)$ or explicitly $\mcR_q(y,y_0) r_q = \mcR_q(x,y_0)$. This is
  equivalent to $r_q \in \mcR_q$. 
\end{proof}

The groupoids defined by different sections may be the same, e. g. if $\mcF$ is
a vector bundle, $\omega$ and $\lambda \omega$ for a constant $\lambda \neq 0$ describe
the same groupoid. Proposition \ref{prop_Groupoid_Bundle} can be extended to
non-homogeneous fibres $F$ with the additional assumption that $\Phi_\omega$ has
constant rank ($\omega(y) 1_y = \omega(y)$ implies $\img(\omega) \subseteq
\img(\Phi_\omega)$). Vessiot \cite{Vessiot03} calls the coordinate expressions of
$\Phi_\omega(r_q) = \omega(s(r_q))$ Lie form. The sequence
(\ref{eq_R_q_sequence}) is due to Pommaret \cite{Pommaret78}.

\section{Systems of PDE and Formal Integrability}
\label{sec_implementation}

The next step is to consider a jet groupoid as a system of PDE (see
e. g. \cite{Goldschmidt67DG}, \cite{KLV86}, \cite{Spencer69}) in order to study
its integrability. The prolongation of a groupoid acting on a fibre bundle was
introduced by Ehresmann \cite{Ehresmann55} (see also \cite{KMS93}).

\begin{definition}
  \label{def_Prolongation_Projection}
  A subbundle $\mcR_q \subseteq J_q(\mcE)$ of the $q$-th order jet
  bundle of a fibre bundle $\mcE \to X$ is called \emph{system of PDE} and solutions are
  (local) sections of $\mcR_q$. The \emph{$r$-prolongation} is the subset
  \[
  \mcR_{q+r} := J_r(\mcR_q) \cap J_{q+r}(\mcE), \qquad r \in \Z_{\geq 0}
  \]
  and 
  \[
  \mcR^{(s)}_{q+r} := \pi^{q+r+s}_{q+r}(\mcR_{q+r+s}) \subseteq \mcR_{q+r}, \qquad r,s \in \Z_{\geq 0}
  \]
  is called \emph{projection}. $\mcR_q$ is called \emph{formally integrable} if
  $\mcR_{q+r}$ is a fibre bundle and the projections $\pi^{q+r+s}_{q+r}:
  \mcR_{q+r+s} \to \mcR_{q+r}$ are surjective submersions for all $r,s \in \Z_{\geq 0}$.
\end{definition}

An effective criterion to decide the formal integrability of a system of PDE was
given by Goldschmidt \cite{Goldschmidt67DG}. It reduces the infinite number of
conditions to a single one, once the \emph{symbol} $g_{q} = (S^q T^*X \otimes
V(\mcE)) \cap V(\mcR_q)$ is $2$-acyclic. For an introduction to symbols and
Spencer cohomology see \cite{Goldschmidt67DG}, \cite{KLV86}, \cite{Spencer69}
or \cite[ch. 7.2]{Pommaret78} for details in the case of jet groupoids. 

\begin{theorem}[{\cite[Thm 8.1]{Goldschmidt67DG}}]
  \label{thm_formalIntCondition}
  Let $\mcR_q \subseteq J_q(\mcE)$ be a system of order $q$ on $\mcE$,
  such that $\mcR_{q+r}$ is a subbundle of $J_{q+r}(\mcE)$. If the symbol $g_q$
  is 2-acyclic, $g_{q+1} \to \mcR_q$ is a vector bundle and if the map $\pi^{q+1}_q :
  \mcR_{q+1} \to \mcR_q$ is surjective, then $\mcR_q$ is formally integrable.
\end{theorem}

We will now derive the bundles and sections that describe the prolongation
$\mcR_{q+r}(\omega)$ and the projection $\mcR^{(1)}_q(\omega)$ of a jet groupoid
$\mcR_q(\omega)$. For a short notation, all maps to the base $X$ (as $s$, $t$ or
$\pi$) will keep their name after prolongating or taking jet bundles. The
prolongation $\mcR_{q+r}(\omega)$ has an obvious description:

\begin{proposition}
  \label{prop_Prolongation_Natural_Bundle}
  Let $\mcF$ be a natural bundle of order $q$ with a section $\omega$. Then
  $J_r(\mcF)$ is a natural bundle of order $q+r$ and the prolongation
  $\mcR_{q+r}(\omega)$ of $\mcR_q(\omega)$ is the symmetry groupoid
  $\Stab_{J_r(\mcF)}^{q+r}(j_r(\omega))$.
\end{proposition}

\begin{proof}
  Apply the functor $J_r()$ to the $\Pi_q$-action on $\mcF$ and use the natural
  embedding $\Pi_{q+r} \hookrightarrow J_r(\Pi_q)$ to establish the
  $\Pi_{q+r}$-action on $J_r(\mcF)$. Note that the image of this embedding is
  $J_r(\Pi_q) \cap \Pi_{q+r}$. As $\mcR_q(\omega)$ is defined as
  $\ker_\omega(\Phi_\omega)$, the exact sequence for $J_r(\mcR_q(\omega))$ is:
  \[
  \xymatrix@C=3pc{
    0 \ar[r] & J_r(\mcR_q(\omega)) \ar[r] & J_r(\Pi_q) \ar@<0.6ex>[r]^{j_r(\Phi_\omega)}
    \ar@<-0.4ex>[r]_{j_r(\omega) \circ s }  &  J_r(\mcF) \\
    0 \ar[r] & \Stab_{J_r(\mcF)}^{q+r}(j_r(\omega)) \ar[r] \ar@{^{(}->}[u] & \Pi_{q+r}
    \ar@<-0.4ex>[r] \ar@<0.6ex>[r] \ar@{^{(}->}[u] & J_r(\mcF) \ar@{=}[u]
  }
  \]
  and the intersection $J_r(\mcR_q(\omega)) \cap \Pi_{q+r}$ actually is the symmetry
  groupoid of $j_r(\omega)$. 
\end{proof}

The fibres of $J_r(\mcF)$ are not necessarily homogeneous, so we cannot assure
that $\mcR_{q+r}(\omega)$ is still a subbundle of $\Pi_{q+r}$ or equivalently a
Lie groupoid. However if the rank of $j_r(\Phi_\omega)$ is constant,
$\mcR_{q+r}(\omega)$ is a Lie groupoid again.

To describe the projections $\mcR^{(1)}_q (\omega)$ we write $J_1(\mcF)$ as a
bundle associated to $P_{q+1}$ with fibre $J_1(F) := J_1(\mcF)_{y_0}$. The idea
to use fibre $F_1 := J_1(F)/K_{q+1}$ to obtain the associated bundle $\mcF_1 := P_{q+1}
\times_{\GL_{q+1}} F_1$ is due to Barakat. 

\begin{proposition}
  \label{prop_Groupoid_Projection}
  $\mcF_1 \cong P_q \times_{\GL_q} F_1$ is a natural bundle of
  order $q$ and if $I : J_1(\mcF) \to \mcF_1$ is the projection,
  $\mcR^{(1)}_q(\omega)$ is the symmetry groupoid $\Stab_{\mcF_1}^q(I(j_1(\omega)))$.
\end{proposition}

\begin{proof}
  By construction of $K_{q+1}$, the $\GL_{q+1}$-action on the fibre $F_1$
  factors over $\GL_q$ and $P_{q+1}/K_{q+1} \cong P_q$ ensures that $\mcF_1$ is
  a natural bundle of order $q$. The preimage $I^{-1}(v)$ of $v \in {\mcF_1}_y$ can be written
  as $u_1 K_{q+1}$ with $u_1 \in I^{-1}(v)$ and $K_{q+1} \unlhd
  \Pi_{q+1}(y,y)$. The action on $\mcF_1$ is defined by $v f_q = u_1 K_{q+1}
  f_{q+1} = u_1 f_{q+1} K_{q+1}$ and we have the exact sequence for
  $\Stab_{\mcF_1}^q(I(j_1(\omega)))$:
  \[
  \xymatrix@C=3pc{
    0 \ar[r] & 
    \Stab_{\mcF_1}^q(I(j_1(\omega)))
    \ar[r] & \Pi_q
    \ar@<-0.4ex>[r]_{I(j_1(\omega)) \circ s} \ar@<0.6ex>[r]^{I(j_1(\Phi_\omega))} & \mcF_1.
  }
  \]
  The symmetry condition
  \[
  I(j_1(\omega))(y) f_q = j_1(\omega)(y) f_{q+1} K_{q+1} \stackrel{!}{=}
  j_1(\omega)(x) K_{q+1}
  \]
  is equivalent to the existence of a preimage $r_{q+1} \in \mcR_{q+1}(\omega)$
  projecting onto $f_q$. 
\end{proof}

We now come to the main result of this article, which is a conceptional proof
using groupoids of a theorem implicitly present in \cite{Vessiot03} and
formulated by Pommaret.

\begin{theorem}
  \label{thm_Projection}
  The projection $\pi^{q+1}_q : \mcR_{q+1}(\omega) \to \mcR_q(\omega)$ is an epimorphism
  if and only if there is a $\Pi_q$-equivariant section $c:\mcF \to \mcF_1$, $c(u f_q) =
  c(u) f_q$, such that $I(j_1(\omega)) = c(\omega)$. This gives the exact
  sequence:
  \[
  \xymatrix@C=3pc{
    0 \ar[r] & 
    \mcR_q(\omega)
    \ar[r] & 
    \Pi_q \ar@<0.6ex>[r]^{\Phi_\omega}  \ar@<-0.4ex>[r]_{\omega \circ s} & 
    \mcF \ar@<0.6ex>[r]^{I \circ j_1}  \ar@<-0.4ex>[r]_{c} & 
    \mcF_1.
  }
  \]
\end{theorem}

\begin{proof}
  Whenever we define an element $a_q$, $a_{q+1}$ denotes an arbitrary preimage
  under the appropriate projection $\pi^{q+1}_q$.
  First assume the existence of $r_{q+1} \in (\pi^{q+1}_q)^{-1} (r_q)$ for all
  $r_q \in \mcR_q(\omega)$. To construct an equivariant section, we define
  \[
  c( \omega(y) ) := j_1(\omega)(y) K_{q+1}.
  \]
  For $\omega(y) \neq u \in F_y$ there is a $g_q \in \GL_q(\RR^n)$ with $u =
  \omega(y) g_q$, we set
  \[
  c( u ) := j_1(\omega)(y) g_{q+1} K_{q+1},
  \]
  which is well-define due to $g_{q+1} K_{q+1}$ being the whole preimage in
  $\GL_{q+1}(\RR^n)$. For each $f_q \in \Pi_q$, we can find $h_q \in \GL_q(\RR^n)$ with
  \[
  \omega (x) h_q = u \, f_q = \omega(y) g_q f_q = \omega(y) r_q h_q \quad
  \text{and} \quad f_q = g_q^{-1} r_q h_q.
  \]
  where the existence of $r_{q+1}$ implies the equivariance:
  \begin{eqnarray*}
    c( u \, f_q ) 
    & = & j_1(\omega)(x) h_{q+1} K_{q+1} \\
    & = & j_1(\omega)(y) g_{q+1} (g_{q+1}^{-1} r_{q+1} h_{q+1} ) K_{q+1} \\
    & = & c(u) f_{q+1} K_{q+1} =  c(u) f_q
  \end{eqnarray*}
  Using the equivariance of $c$ on $c(\omega(y) \, r_q) = c( \omega(y) ) r_q$,
  we obtain
  \[
  j_1(\omega)(y) \bar{r}_{q+1} K_{q+1} = j_1(\omega)(x) K_{q+1}
  \]
  for an arbitrary preimage $\bar{r}_{q+1}$. There is a $k_{q+1}$ such that
  $r_{q+1} = \bar{r}_{q+1} k_{q+1}$ satisfying
  \[
  j_1(\omega)(y) r_{q+1} = j_1(\omega)(x)
  \]
  which provides a preimage $r_{q+1} \in \mcR_{q+1}(\omega)$ for $r_q$. 
\end{proof}


Using the $\GL_q$-action on the fibre $F_1$, all possibilities for
equivariant sections $c$ can be calculated. The resulting integrability
conditions $I(j_1(\omega)) = c (\omega)$ are called \emph{Vessiot structure
  equations}. They express the condition that each defining equation for
$\mcR_{q+1}(\omega)$ where the jets of order $q+1$ can be eliminated must be a
consequence of the equations for $\mcR_q(\omega)$.

If the Vessiot structure equations are fulfilled for a section $\omega$, 
$\mcR_{q+1}(\omega)$ is transitive and a subbundle of $\Pi_{q+1}$. Then by
\cite{Pommaret78}, the symbol $g_{q+1}$ is a vector bundle and we can apply
Theorem \ref{thm_formalIntCondition} which implies formal integrability. Theorem
\ref{thm_Projection} can be extended to non-homogeneous fibres $F$ as long as
the section $\omega$ defines a Lie groupoid.

Starting from an arbitrary transitive jet groupoid $\mcR_q$, we have found a
natural bundle $\mcF$ of geometric objects and a special object $\omega_0$,
such that $\mcR_q = \Stab_\mcF^q(\omega_0)$. 
It has been shown that the section $j_r(\omega)$ of $J_r(\mcF)$ defines the
$r$-th prolongation $\mcR_{q+r}(\omega)$ and that $I (j_1 (\omega))$ on $\mcF_1$
corresponds to the projection $\mcR^{(1)}_q(\omega)$. Based on Theorem
\ref{thm_formalIntCondition}, the projection theorem leads to a check of formal
integrability directly on the level of sections $\omega$ of $\mcF$. In most
cases, the integrability conditions have an immediate geometric interpretation
as in the following example.

\section{Example}
\label{sec_example}

\begin{maplegroup}
The following calculation is due to Barakat using the {\sc Maple}
package {\tt jets} \cite{Barakat01}, which contains routines for jet
groupoids and natural bundles. It will be used to show explicit
examples of the objects in the theoretical part.

\begin{mapleinput}
\mapleinline{active}{1d}{with(jets):}{%
}
\end{mapleinput}

\end{maplegroup}
\begin{maplegroup}
\noindent Dimension of the base manifold $X$ and some coordinates:

\end{maplegroup}
\begin{maplegroup}
\begin{mapleinput}
\mapleinline{active}{1d}{n:=2: 
ivar:=[x1,x2]: dvar:=[y1,y2]: 
Ivar:=[phi1,phi2]: Dvar:=[xi1,xi2]:}{%
}
\end{mapleinput}

\end{maplegroup}
\begin{maplegroup}
\noindent The jet groupoid expressing the invariance of the flat
euclidean metric $g$ on $X$:

\end{maplegroup}
\begin{maplegroup}
\begin{mapleinput}
\mapleinline{active}{1d}{(Jac,g):=(matrix(n,n,jetcoor(1,ivar,dvar)), linalg[diag](1$n));}{%
}
\end{mapleinput}

\mapleresult
\begin{maplelatex}
\mapleinline{inert}{2d}{Jac, g := matrix([[y1[x1], y1[x2]], [y2[x1], y2[x2]]]), matrix([[1,
0], [0, 1]]);}{%
\[
\mathit{Jac}, \,g :=  \left[ 
{\begin{array}{cc}
{\mathit{y1}_{\mathit{x1}}} & {\mathit{y1}_{\mathit{x2}}} \\
{\mathit{y2}_{\mathit{x1}}} & {\mathit{y2}_{\mathit{x2}}}
\end{array}}
 \right] , \, \left[ 
{\begin{array}{rr}
1 & 0 \\
0 & 1
\end{array}}
 \right] 
\]
}
\end{maplelatex}

\end{maplegroup}
\begin{maplegroup}
\begin{mapleinput}
\mapleinline{active}{1d}{J_:=evalm(linalg[transpose](Jac) &* g &* Jac);}{%
}
\end{mapleinput}

\mapleresult
\begin{maplelatex}
\mapleinline{inert}{2d}{J_ := matrix([[y1[x1]^2+y2[x1]^2, y1[x1]*y1[x2]+y2[x1]*y2[x2]],
[y1[x1]*y1[x2]+y2[x1]*y2[x2], y1[x2]^2+y2[x2]^2]]);}{%
\[
\mathit{J\_} :=  \left[ 
{\begin{array}{cc}
{\mathit{y1}_{\mathit{x1}}}^{2} + {\mathit{y2}_{\mathit{x1}}}^{2}
 & {\mathit{y1}_{\mathit{x1}}}\,{\mathit{y1}_{\mathit{x2}}} + {
\mathit{y2}_{\mathit{x1}}}\,{\mathit{y2}_{\mathit{x2}}} \\
{\mathit{y1}_{\mathit{x1}}}\,{\mathit{y1}_{\mathit{x2}}} + {
\mathit{y2}_{\mathit{x1}}}\,{\mathit{y2}_{\mathit{x2}}} & {
\mathit{y1}_{\mathit{x2}}}^{2} + {\mathit{y2}_{\mathit{x2}}}^{2}
\end{array}}
 \right] 
\]
}
\end{maplelatex}

\end{maplegroup}
\begin{maplegroup}
\begin{mapleinput}
\mapleinline{active}{1d}{GR_g:=[ J_[1,1]=1, J_[1,2]=0, J_[2,2]=1];}{%
}
\end{mapleinput}

\mapleresult
\begin{maplelatex}
\mapleinline{inert}{2d}{GR_g := [y1[x1]^2+y2[x1]^2 = 1, y1[x1]*y1[x2]+y2[x1]*y2[x2] = 0,
y1[x2]^2+y2[x2]^2 = 1];}{%
\[
\mathit{GR\_g} := [{\mathit{y1}_{\mathit{x1}}}^{2} + {\mathit{y2}
_{\mathit{x1}}}^{2}=1, \,{\mathit{y1}_{\mathit{x1}}}\,{\mathit{y1
}_{\mathit{x2}}} + {\mathit{y2}_{\mathit{x1}}}\,{\mathit{y2}_{
\mathit{x2}}}=0, \,{\mathit{y1}_{\mathit{x2}}}^{2} + {\mathit{y2}
_{\mathit{x2}}}^{2}=1]
\]
}
\end{maplelatex}

\end{maplegroup}
\begin{maplegroup}
These equations locally define a transitive groupoid $\mathcal{R}_1(g)
\subset \Pi_1$ with isotropy groups $O_2(\mathbb{R})$. They have been
constructed by the action of $\GL_1 \cong \GL(\mathbb{R}^2)$ on the
space $F$ of scalar products on $\mathbb{R}^2$. So we start with the 
natural bundle $\mathcal{F}_g  = P_1 \times_{\GL_1} F \cong S^2 T^*
X_{\geq 0}$ of symmetric positive definite 2-forms and the equations
for $\mathcal{R}_1(g)$ are already in Lie form (see section
\ref{sec_Jet_Groupoids} for $P_q$ and $\GL_q$). Define coordinates for
$\mathcal{F}_g$ and a section $\omega$: 

\end{maplegroup}
\begin{maplegroup}
\begin{mapleinput}
\mapleinline{active}{1d}{uvar_g:=[u11,u12,u22]: wvar_g:=[omega11,omega12,omega22]:}{%
}
\end{mapleinput}

\end{maplegroup}
\begin{maplegroup}
\noindent As in \cite{Vessiot03}, the coordinate changes of
$\mathcal{F}_g$ are given in the form \[ \hat{x} = \phi(x), \quad u =
\Psi(\hat{x}=\phi(x), \hat{u}, \phi_q(x)) \] (mind the hats in the
second equation). For shorter output, jet notation is used for
$\phi(x)$ and its derivatives:

\end{maplegroup}
\begin{maplegroup}
\begin{mapleinput}
\mapleinline{active}{1d}{inv_g:=ezip(uvar_g,map(lhs,GR_g)):
F_g:=natfin(inv_g,ivar,dvar,uvar_g,Ivar,""):
eqn2ind(F_g,ivar,Ivar);}{%
}
\end{mapleinput}

\mapleresult
\begin{maplelatex}
\mapleinline{inert}{2d}{[x1 = phi1, x2 = phi2, u11 =
phi1[x1]^2*u11+2*phi1[x1]*phi2[x1]*u12+phi2[x1]^2*u22, u12 =
phi1[x2]*phi1[x1]*u11+phi1[x2]*phi2[x1]*u12+phi2[x2]*phi1[x1]*u12+phi2
[x2]*phi2[x1]*u22, u22 =
phi1[x2]^2*u11+2*phi1[x2]*phi2[x2]*u12+phi2[x2]^2*u22];}{%
\maplemultiline{
[\mathit{x1}=\phi 1, \,\mathit{x2}=\phi 2, \\
\mathit{u11}={\phi 1_{\mathit{x1}}}^{2}\,\mathit{u11} + 2\,{\phi
  1_{\mathit{x1}}}\,{\phi 2_{\mathit{x1}}}\,\mathit{u12} + {\phi
  2_{\mathit{x1}}}^{2} \,\mathit{u22},  \\
\mathit{u12}={\phi 1_{\mathit{x2}}}\,{\phi 1_{\mathit{x1}}}\,
\mathit{u11} + {\phi 1_{\mathit{x2}}}\,{\phi 2_{\mathit{x1}}}\,
\mathit{u12} + {\phi 2_{\mathit{x2}}}\,{\phi 1_{\mathit{x1}}}\,
\mathit{u12} + {\phi 2_{\mathit{x2}}}\,{\phi 2_{\mathit{x1}}}\,
\mathit{u22},  \\
\mathit{u22}={\phi 1_{\mathit{x2}}}^{2}\,\mathit{u11} + 2\,{\phi 
1_{\mathit{x2}}}\,{\phi 2_{\mathit{x2}}}\,\mathit{u12} + {\phi 2
_{\mathit{x2}}}^{2}\,\mathit{u22}] }
}
\end{maplelatex}

\end{maplegroup}
\begin{maplegroup}
\noindent The groupoid $\mathcal{R}_1 (\omega)$ for a general section
in Lie form:

\end{maplegroup}
\begin{maplegroup}
\begin{mapleinput}
\mapleinline{active}{1d}{LieFormG(F_g,ivar,dvar,Ivar,wvar_g);}{%
}
\end{mapleinput}

\mapleresult
\begin{maplelatex}
\mapleinline{inert}{2d}{[y1[x1]^2*omega11(y1,y2)+2*y1[x1]*y2[x1]*omega12(y1,y2)+y2[x1]^2*omeg
a22(y1,y2) = omega11(x1,x2),
y1[x2]*y1[x1]*omega11(y1,y2)+y1[x2]*y2[x1]*omega12(y1,y2)+y2[x2]*y1[x1
]*omega12(y1,y2)+y2[x2]*y2[x1]*omega22(y1,y2) = omega12(x1,x2),
y1[x2]^2*omega11(y1,y2)+2*y1[x2]*y2[x2]*omega12(y1,y2)+y2[x2]^2*omega2
2(y1,y2) = omega22(x1,x2)];}{%
\maplemultiline{
[{\mathit{y1}_{\mathit{x1}}}^{2}\,\omega 11(\mathit{y1}, \,
\mathit{y2}) + 2\,{\mathit{y1}_{\mathit{x1}}}\,{\mathit{y2}_{
\mathit{x1}}}\,\omega 12(\mathit{y1}, \,\mathit{y2}) + {\mathit{
y2}_{\mathit{x1}}}^{2}\,\omega 22(\mathit{y1}, \,\mathit{y2})=
\omega 11(\mathit{x1}, \,\mathit{x2}),  \\
{\mathit{y1}_{\mathit{x2}}}\,{\mathit{y1}_{\mathit{x1}}}\,\omega 
11(\mathit{y1}, \,\mathit{y2}) + {\mathit{y1}_{\mathit{x2}}}\,{
\mathit{y2}_{\mathit{x1}}}\,\omega 12(\mathit{y1}, \,\mathit{y2})
 + {\mathit{y2}_{\mathit{x2}}}\,{\mathit{y1}_{\mathit{x1}}}\,
\omega 12(\mathit{y1}, \,\mathit{y2}) \\
\mbox{} + {\mathit{y2}_{\mathit{x2}}}\,{\mathit{y2}_{\mathit{x1}}
}\,\omega 22(\mathit{y1}, \,\mathit{y2})=\omega 12(\mathit{x1}, 
\,\mathit{x2}),  \\
{\mathit{y1}_{\mathit{x2}}}^{2}\,\omega 11(\mathit{y1}, \,
\mathit{y2}) + 2\,{\mathit{y1}_{\mathit{x2}}}\,{\mathit{y2}_{
\mathit{x2}}}\,\omega 12(\mathit{y1}, \,\mathit{y2}) + {\mathit{
y2}_{\mathit{x2}}}^{2}\,\omega 22(\mathit{y1}, \,\mathit{y2})=
\omega 22(\mathit{x1}, \,\mathit{x2})] }
}
\end{maplelatex}

\end{maplegroup}
\begin{maplegroup}
\noindent The special section $\omega_0$ for the flat metric $g$:

\begin{mapleinput}
\mapleinline{active}{1d}{omega0:=map(rhs,GR_g);}{%
}
\end{mapleinput}

\mapleresult
\begin{maplelatex}
\mapleinline{inert}{2d}{omega0 := [1, 0, 1];}{%
\[
\omega 0 := [1, \,0, \,1]
\]
}
\end{maplelatex}

\end{maplegroup}
\begin{maplegroup}
The application of Theorem \ref{thm_Projection} at this point gives no
integrability conditions, although an arbitrary metric should not be
integrable. The reason is that the symbol of $\mathcal{R}_1(\omega)$
is not yet $2$-acyclic, but $\mathcal{R}_2(\omega)$ has $2$-acyclic
symbol. We could go on with $J_1(\mathcal{F}_g)$, but in order to keep
geometrical interpretation (and short expressions) we also model the
Christoffel symbols by plugging the derivatives of the transformed 
flat metric (\texttt{GR\_g}) into: \begin{equation}
\label{eq_Christoffel_Symbol} \Gamma^k_{ij} (x) = \frac{1}{2} g^{k
\mu} (x) \left( \frac{\partial g_{i\mu}}{\partial x^j}(x) +
        \frac{\partial g_{j\mu}}{\partial x^i}(x) - \frac{\partial
g_{ij}}{\partial x^\mu}(x) \right) ,  \end{equation} which gives the
equations for the Christoffel symbols of the flat metric in Lie form:

\end{maplegroup}
\begin{maplegroup}
\begin{mapleinput}
\mapleinline{active}{1d}{dJac := linalg[det](Jac):
Phi_Gamma := [
(y2[x2]*y1[x1,x1]-y2[x1,x1]*y1[x2])/dJac,
(y2[x2]*y1[x1,x2]-y2[x1,x2]*y1[x2])/dJac,
(y2[x1,x1]*y1[x1]-y1[x1,x1]*y2[x1])/dJac,
(y2[x1,x2]*y1[x1]-y1[x1,x2]*y2[x1])/dJac,
(y2[x2]*y1[x2,x2]-y2[x2,x2]*y1[x2])/dJac,
(y2[x2,x2]*y1[x1]-y1[x2,x2]*y2[x1])/dJac]:}{%
}
\end{mapleinput}

\end{maplegroup}
\begin{maplegroup}
\noindent The coordinates for the Christoffel symbols ({\texttt uijk}
stands for $\Gamma^i_{jk}$):

\end{maplegroup}
\begin{maplegroup}
\begin{mapleinput}
\mapleinline{active}{1d}{uvar_Gamma:=[u111,u112,u211,u212,u122,u222]:}{%
}
\end{mapleinput}

\end{maplegroup}
\begin{maplegroup}
\begin{mapleinput}
\end{mapleinput}

\end{maplegroup}
\begin{maplegroup}
\noindent Calculate the natural bundle $\mathcal{F}_\Gamma$ of
Christoffel symbols:

\begin{mapleinput}
\mapleinline{active}{1d}{inv_Gamma := ezip(uvar_Gamma,Phi_Gamma):
F_Gamma:=natfin(inv_Gamma,ivar,dvar,uvar_Gamma,Ivar,""):}{%
}
\end{mapleinput}

\end{maplegroup}
\begin{maplegroup}
\noindent The result is $\mathcal{F} =  \mathcal{F}_g \times_X
\mathcal{F}_\Gamma \cong J_1(\mathcal{F}_g)$. Usually,
$J_1(\mathcal{F})$ is only an affine bundle over $\mathcal{F}$ and
does not split. The fibre $F$ is a homogeneous $\GL_2$-space, so each
section on $\mathcal{F}$ defines a Lie groupoid.

\begin{mapleinput}
\mapleinline{active}{1d}{uvar:=[op(uvar_g),op(uvar_Gamma)]: 
F:=[op(F_g),op(F_Gamma[n+1..-1])]:}{%
}
\end{mapleinput}

\end{maplegroup}
\begin{maplegroup}
To calculate the projection to the bundle $\mathcal{F}_1$,  the vector
fields of infinitesimal transformations of $\mathcal{F}$. If $\xi^i(x)
\frac{\partial}{\partial x^i}$ is a vector field on $X$, it can be
extended to $\mathcal{F}$:

\end{maplegroup}
\begin{maplegroup}
\begin{mapleinput}
\mapleinline{active}{1d}{vec:=natfin2inf(F,ivar,Ivar,Dvar,"");}{%
}
\end{mapleinput}

\mapleresult
\begin{maplelatex}
\mapleinline{inert}{2d}{vec := [[xi1, [x1]], [xi2, [x2]], [-2*u11*xi1[x1]-2*u12*xi2[x1],
[u11]], [-xi1[x1]*u12-xi1[x2]*u11-xi2[x1]*u22-xi2[x2]*u12, [u12]],
[-2*u12*xi1[x2]-2*u22*xi2[x2], [u22]],
[-xi1[x1,x1]-xi1[x1]*u111+xi1[x2]*u211-2*xi2[x1]*u112, [u111]],
[-xi1[x1,x2]-xi1[x2]*u111+xi1[x2]*u212-xi2[x1]*u122-xi2[x2]*u112,
[u112]],
[-2*xi1[x1]*u211-xi2[x1,x1]+xi2[x1]*u111-2*xi2[x1]*u212+xi2[x2]*u211,
[u211]],
[-xi1[x1]*u212-xi1[x2]*u211-xi2[x1,x2]+xi2[x1]*u112-xi2[x1]*u222,
[u212]],
[-xi1[x2,x2]+xi1[x1]*u122-2*xi1[x2]*u112+xi1[x2]*u222-2*xi2[x2]*u122,
[u122]], [-2*xi1[x2]*u212-xi2[x2,x2]+xi2[x1]*u122-xi2[x2]*u222,
[u222]]];}{%
\maplemultiline{
\mathit{vec} := [[\xi 1, \,[\mathit{x1}]], \,[\xi 2, \,[\mathit{
x2}]], \,[ - 2\,\mathit{u11}\,{\xi 1_{\mathit{x1}}} - 2\,\mathit{
u12}\,{\xi 2_{\mathit{x1}}}, \,[\mathit{u11}]],  \\
[ - {\xi 1_{\mathit{x1}}}\,\mathit{u12} - {\xi 1_{\mathit{x2}}}\,
\mathit{u11} - {\xi 2_{\mathit{x1}}}\,\mathit{u22} - {\xi 2_{
\mathit{x2}}}\,\mathit{u12}, \,[\mathit{u12}]],  \\
[ - 2\,\mathit{u12}\,{\xi 1_{\mathit{x2}}} - 2\,\mathit{u22}\,{
\xi 2_{\mathit{x2}}}, \,[\mathit{u22}]],  \\
[ - {\xi 1_{\mathit{x1}, \,\mathit{x1}}} - {\xi 1_{\mathit{x1}}}
\,\mathit{u111} + {\xi 1_{\mathit{x2}}}\,\mathit{u211} - 2\,{\xi 
2_{\mathit{x1}}}\,\mathit{u112}, \,[\mathit{u111}]],  \\
[ - {\xi 1_{\mathit{x1}, \,\mathit{x2}}} - {\xi 1_{\mathit{x2}}}
\,\mathit{u111} + {\xi 1_{\mathit{x2}}}\,\mathit{u212} - {\xi 2_{
\mathit{x1}}}\,\mathit{u122} - {\xi 2_{\mathit{x2}}}\,\mathit{
u112}, \,[\mathit{u112}]],  \\
[ - 2\,{\xi 1_{\mathit{x1}}}\,\mathit{u211} - {\xi 2_{\mathit{x1}
, \,\mathit{x1}}} + {\xi 2_{\mathit{x1}}}\,\mathit{u111} - 2\,{
\xi 2_{\mathit{x1}}}\,\mathit{u212} + {\xi 2_{\mathit{x2}}}\,
\mathit{u211}, \,[\mathit{u211}]],  \\
[ - {\xi 1_{\mathit{x1}}}\,\mathit{u212} - {\xi 1_{\mathit{x2}}}
\,\mathit{u211} - {\xi 2_{\mathit{x1}, \,\mathit{x2}}} + {\xi 2_{
\mathit{x1}}}\,\mathit{u112} - {\xi 2_{\mathit{x1}}}\,\mathit{
u222}, \,[\mathit{u212}]],  \\
[ - {\xi 1_{\mathit{x2}, \,\mathit{x2}}} + {\xi 1_{\mathit{x1}}}
\,\mathit{u122} - 2\,{\xi 1_{\mathit{x2}}}\,\mathit{u112} + {\xi 
1_{\mathit{x2}}}\,\mathit{u222} - 2\,{\xi 2_{\mathit{x2}}}\,
\mathit{u122}, \,[\mathit{u122}]],  \\
[ - 2\,{\xi 1_{\mathit{x2}}}\,\mathit{u212} - {\xi 2_{\mathit{x2}
, \,\mathit{x2}}} + {\xi 2_{\mathit{x1}}}\,\mathit{u122} - {\xi 2
_{\mathit{x2}}}\,\mathit{u222}, \,[\mathit{u222}]]] }
}
\end{maplelatex}

\end{maplegroup}
\begin{maplegroup}
The above list denotes a vector field. It is read as follows: $[\xi 1,
[x1]]$ stands for $\xi^1(x) \frac{\partial}{\partial x^1}$ and the
complete vector field is obtained by adding up all list entries. Each
choice of $\xi^i, \ldots, \xi^i_{x^i,x^j}$ gives an infinitesimal
transformation of $\mathcal{F}$. We calculate the coordinates of
$\mathcal{F}_1$ that express the projection $I: J_1(\mathcal{F}) \to \mathcal{F}_1$:

\end{maplegroup}
\begin{maplegroup}
\begin{mapleinput}
\mapleinline{active}{1d}{F1:=F1coor(vec,ivar,Dvar,uvar);}{%
}
\end{mapleinput}

\mapleresult
\begin{maplelatex}
\mapleinline{inert}{2d}{F1 := [u11[x1], u11[x2], u12[x1], u12[x2], u22[x1], u22[x2],
u111[x2]-u112[x1], u112[x2]-u122[x1], u211[x2]-u212[x1],
u212[x2]-u222[x1]];}{%
\maplemultiline{
\mathit{F1} := [{\mathit{u11}_{\mathit{x1}}}, \,{\mathit{u11}_{
\mathit{x2}}}, \,{\mathit{u12}_{\mathit{x1}}}, \,{\mathit{u12}_{
\mathit{x2}}}, \,{\mathit{u22}_{\mathit{x1}}}, \,{\mathit{u22}_{
\mathit{x2}}}, \\ 
\qquad {\mathit{u111}_{\mathit{x2}}} - {\mathit{u112}_{
\mathit{x1}}}, \,{\mathit{u112}_{\mathit{x2}}} - {\mathit{u122}_{
\mathit{x1}}}, {\mathit{u211}_{\mathit{x2}}} - {\mathit{u212}_{\mathit{x1}}}, \,
{\mathit{u212}_{\mathit{x2}}} - {\mathit{u222}_{\mathit{x1}}}] }
}
\end{maplelatex}

\end{maplegroup}
\begin{maplegroup}
\begin{mapleinput}
\mapleinline{active}{1d}{d1:=nops(F1):
vvar := [v1,v2,v3,v4,v5,v6,v7,v8,v9,v10]:}{%
}
\end{mapleinput}

\end{maplegroup}
\begin{maplegroup}
\noindent All further computations only need the infinitesimal
coordinate changes of $\mathcal{F}_1$. Setting zero order jets $\xi^i
= 0$ to zero and collecting for higher order jets of $\xi^i$, the list
\texttt{L1} contains a representation of the Lie algebra of
$\GL_2(\RR^2)$ as vertical vector fields on $\mathcal{F}_1$. 

\end{maplegroup}
\begin{maplegroup}
\begin{mapleinput}
\mapleinline{active}{1d}{inv1 := ezip(vvar,F1):
vec1:=natinfG(vec,inv1,ivar,uvar,vvar,Dvar):
L1:=lstvec(sortcon(vec1,[op(uvar),op(vvar)]),ivar,Dvar,""):}{%
}
\end{mapleinput}

\end{maplegroup}
\begin{maplegroup}
Before calculating the possible equivariant sections of
$\mathcal{F}_1$, we will modify the coordinates of $\mathcal{F}_1$ to
obtain a vector bundle atlas. This is achieved by choosing the coordinates for
the fibres of $\mathcal{F}_1 \to \mathcal{F}$ to be $K_2$-invariant (see
sequence (\ref{eq_GL_exact_seq}) for $K_{q+1}$):

\end{maplegroup}
\begin{maplegroup}
\begin{mapleinput}
\mapleinline{active}{1d}{cvar := [c1,c2,c3,c4,c5,c6,c7,c8,c9,c10]:
subv:=map(a->vvar[a]=cvar[a](op(uvar)),[$1..nops(cvar)]):
cndi:=map(i->subs(subv,invcond(L1[1][n^2+1..-1],
                  [lhs(subv[i])-rhs(subv[i])],L1[2])[1]),[$1..d1]):
sol_Gamma := map(ci->jsolve(cndi[ci], uvar,
                  [cvar[ci](op(uvar))],""), [$1..d1]):}{%
}
\end{mapleinput}

\end{maplegroup}
\begin{maplegroup}
\begin{mapleinput}
\end{mapleinput}

\end{maplegroup}
\begin{maplegroup}
\noindent The results all depend on arbitrary functions $\_F1(u11,
u12, u22)$, which will be set to zero:

\end{maplegroup}
\begin{maplegroup}
\begin{mapleinput}
\mapleinline{active}{1d}{sol_Gamma[1];}{%
}
\end{mapleinput}

\mapleresult
\begin{maplelatex}
\mapleinline{inert}{2d}{[c1(u11,u12,u22,u111,u112,u211,u212,u122,u222) =
2*u11*u111+2*u12*u211+_F1(u11,u12,u22)];}{%
\maplemultiline{
[\mathrm{c1}(\mathit{u11}, \,\mathit{u12}, \,\mathit{u22}, \,
\mathit{u111}, \,\mathit{u112}, \,\mathit{u211}, \,\mathit{u212}
, \,\mathit{u122}, \,\mathit{u222})= \\
2\,\mathit{u11}\,\mathit{u111} + 2\,\mathit{u12}\,\mathit{u211}
 + \mathrm{\_F1}(\mathit{u11}, \,\mathit{u12}, \,\mathit{u22})]
 }
}
\end{maplelatex}

\end{maplegroup}
\begin{maplegroup}
\begin{mapleinput}
\mapleinline{active}{1d}{sol_Gamma := eval(map(a->op(subs(_F1=0,a)),sol_Gamma)):}{%
}
\end{mapleinput}

\end{maplegroup}
\begin{maplegroup}
\begin{mapleinput}
\end{mapleinput}

\end{maplegroup}
\begin{maplegroup}
\noindent The new infinitesimal coordinate changes show the vector
bundle structure of $\mathcal{F}_1 \to \mathcal{F}$:

\end{maplegroup}
\begin{maplegroup}
\begin{mapleinput}
\mapleinline{active}{1d}{inv1_1 := zip((a,b)->lhs(a) = rhs(a) - rhs(b),inv1,sol_Gamma):
vec1_1:=natinfG(vec,inv1_1,ivar,uvar,vvar,Dvar);
L1_1:=lstvec(sortcon(vec1_1,[op(uvar),op(vvar)]),ivar,Dvar,""):}{%
}
\end{mapleinput}

\mapleresult
\begin{maplelatex}
\mapleinline{inert}{2d}{vec1_1 := [[xi1, [x1]], [xi2, [x2]], [-2*u11*xi1[x1]-2*u12*xi2[x1],
[u11]], [-xi1[x1]*u12-xi1[x2]*u11-xi2[x1]*u22-xi2[x2]*u12, [u12]],
[-2*u12*xi1[x2]-2*u22*xi2[x2], [u22]],
[-xi1[x1,x1]-xi1[x1]*u111+xi1[x2]*u211-2*xi2[x1]*u112, [u111]],
[-xi1[x1,x2]-xi1[x2]*u111+xi1[x2]*u212-xi2[x1]*u122-xi2[x2]*u112,
[u112]],
[-2*xi1[x1]*u211-xi2[x1,x1]+xi2[x1]*u111-2*xi2[x1]*u212+xi2[x2]*u211,
[u211]],
[-xi1[x1]*u212-xi1[x2]*u211-xi2[x1,x2]+xi2[x1]*u112-xi2[x1]*u222,
[u212]],
[-xi1[x2,x2]+xi1[x1]*u122-2*xi1[x2]*u112+xi1[x2]*u222-2*xi2[x2]*u122,
[u122]], [-2*xi1[x2]*u212-xi2[x2,x2]+xi2[x1]*u122-xi2[x2]*u222,
[u222]], [-3*v1*xi1[x1]-v2*xi2[x1]-2*v3*xi2[x1], [v1]],
[-v1*xi1[x2]-v2*xi2[x2]-2*v2*xi1[x1]-2*v4*xi2[x1], [v2]],
[-v1*xi1[x2]-xi2[x2]*v3-2*v3*xi1[x1]-v4*xi2[x1]-xi2[x1]*v5, [v3]],
[-xi1[x2]*v2-v3*xi1[x2]-2*v4*xi2[x2]-xi1[x1]*v4-v6*xi2[x1], [v4]],
[-2*v3*xi1[x2]-2*v5*xi2[x2]-v5*xi1[x1]-v6*xi2[x1], [v5]],
[-2*v4*xi1[x2]-v5*xi1[x2]-3*v6*xi2[x2], [v6]],
[-xi2[x2]*v7-xi1[x1]*v7-xi2[x1]*v8+xi1[x2]*v9, [v7]],
[-xi1[x2]*v7-2*xi2[x2]*v8+xi1[x2]*v10, [v8]],
[xi2[x1]*v7-2*xi1[x1]*v9-xi2[x1]*v10, [v9]],
[xi2[x1]*v8-xi1[x2]*v9-xi2[x2]*v10-xi1[x1]*v10, [v10]]];}{%
\maplemultiline{
\mathit{vec1\_1} := [
\ldots \mathit{vec} \ldots \\
[ - 3\,\mathit{v1}\,{\xi 1_{\mathit{x1}}} - \mathit{v2}\,{\xi 2_{
\mathit{x1}}} - 2\,\mathit{v3}\,{\xi 2_{\mathit{x1}}}, \,[
\mathit{v1}]],  \\
[ - \mathit{v1}\,{\xi 1_{\mathit{x2}}} - \mathit{v2}\,{\xi 2_{
\mathit{x2}}} - 2\,\mathit{v2}\,{\xi 1_{\mathit{x1}}} - 2\,
\mathit{v4}\,{\xi 2_{\mathit{x1}}}, \,[\mathit{v2}]],  \\
[ - \mathit{v1}\,{\xi 1_{\mathit{x2}}} - {\xi 2_{\mathit{x2}}}\,
\mathit{v3} - 2\,\mathit{v3}\,{\xi 1_{\mathit{x1}}} - \mathit{v4}
\,{\xi 2_{\mathit{x1}}} - {\xi 2_{\mathit{x1}}}\,\mathit{v5}, \,[
\mathit{v3}]],  \\
[ - {\xi 1_{\mathit{x2}}}\,\mathit{v2} - \mathit{v3}\,{\xi 1_{
\mathit{x2}}} - 2\,\mathit{v4}\,{\xi 2_{\mathit{x2}}} - {\xi 1_{
\mathit{x1}}}\,\mathit{v4} - \mathit{v6}\,{\xi 2_{\mathit{x1}}}, 
\,[\mathit{v4}]],  \\
[ - 2\,\mathit{v3}\,{\xi 1_{\mathit{x2}}} - 2\,\mathit{v5}\,{\xi 
2_{\mathit{x2}}} - \mathit{v5}\,{\xi 1_{\mathit{x1}}} - \mathit{
v6}\,{\xi 2_{\mathit{x1}}}, \,[\mathit{v5}]],  \\
[ - 2\,\mathit{v4}\,{\xi 1_{\mathit{x2}}} - \mathit{v5}\,{\xi 1_{
\mathit{x2}}} - 3\,\mathit{v6}\,{\xi 2_{\mathit{x2}}}, \,[
\mathit{v6}]],  \\
[ - {\xi 2_{\mathit{x2}}}\,\mathit{v7} - {\xi 1_{\mathit{x1}}}\,
\mathit{v7} - {\xi 2_{\mathit{x1}}}\,\mathit{v8} + {\xi 1_{
\mathit{x2}}}\,\mathit{v9}, \,[\mathit{v7}]],  \\
[ - {\xi 1_{\mathit{x2}}}\,\mathit{v7} - 2\,{\xi 2_{\mathit{x2}}}
\,\mathit{v8} + {\xi 1_{\mathit{x2}}}\,\mathit{v10}, \,[\mathit{
v8}]], \\
[{\xi 2_{\mathit{x1}}}\,\mathit{v7} - 2\,{\xi 1_{
\mathit{x1}}}\,\mathit{v9} - {\xi 2_{\mathit{x1}}}\,\mathit{v10}
, \,[\mathit{v9}]],  \\
[{\xi 2_{\mathit{x1}}}\,\mathit{v8} - {\xi 1_{\mathit{x2}}}\,
\mathit{v9} - {\xi 2_{\mathit{x2}}}\,\mathit{v10} - {\xi 1_{
\mathit{x1}}}\,\mathit{v10}, \,[\mathit{v10}]]] }
}
\end{maplelatex}

\end{maplegroup}
\begin{maplegroup}
\noindent We are now able to calculate all equivariant sections $c$.
The infinitesimal conditions for equivariance are obtained by applying
all vector fields of the Lie algebra of $\GL_2$ to $v_i - c_i(u) = 0$
and then substituting $v_i \to c_i$. The same method was used for
the $K_2$-invariance.

\end{maplegroup}
\begin{maplegroup}
\begin{mapleinput}
\mapleinline{active}{1d}{subv_1:=map(a->vvar[a]=cvar[a](op(uvar)),[$1..nops(cvar)]):
cnd1 := subs(subv_1,invcond(L1_1[1],map(a->lhs(a)-rhs(a),subv_1),
             L1_1[2])[1]):
cc := jsolve(cnd1,uvar,map(a->a(op(uvar)),cvar),""):}{%
}
\end{mapleinput}

\end{maplegroup}
\begin{maplegroup}
\noindent The Vessiot structure equations show the integrability
conditions with the equivariant sections on the right hand side:

\end{maplegroup}
\begin{maplegroup}
\begin{mapleinput}
\mapleinline{active}{1d}{Ves := subs(inv1_1,subs(cc,subv_1));}{%
}
\end{mapleinput}

\mapleresult
\begin{maplelatex}
\mapleinline{inert}{2d}{Ves := [u11[x1]-2*u11*u111-2*u12*u211 = 0,
u11[x2]-2*u11*u112-2*u12*u212 = 0,
u12[x1]-(u111+u212)*u12-u11*u112-u22*u211 = 0,
u12[x2]-(u112+u222)*u12-u22*u212-u11*u122 = 0,
u22[x1]-2*u12*u112-2*u22*u212 = 0, u22[x2]-2*u12*u122-2*u22*u222 = 0,
u111[x2]-u112[x1]-u212*u112+u122*u211 =
_C1*u12+(-u12^2+u22*u11)^(1/2)*_C2,
u112[x2]-u122[x1]-(u111-u212)*u122-u112*u222+u112^2 = _C1*u22,
u211[x2]-u212[x1]-u212^2+u212*u111-(u112-u222)*u211 = -u11*_C1,
u212[x2]-u222[x1]-u122*u211+u212*u112 =
-_C1*u12+(-u12^2+u22*u11)^(1/2)*_C2];}{%
\maplemultiline{
\mathit{Ves} := [{\mathit{u11}_{\mathit{x1}}} - 2\,\mathit{u11}\,
\mathit{u111} - 2\,\mathit{u12}\,\mathit{u211}=0, \\
{\mathit{u11}
_{\mathit{x2}}} - 2\,\mathit{u11}\,\mathit{u112} - 2\,\mathit{u12
}\,\mathit{u212}=0,  \\
{\mathit{u12}_{\mathit{x1}}} - (\mathit{u111} + \mathit{u212})\,
\mathit{u12} - \mathit{u11}\,\mathit{u112} - \mathit{u22}\,
\mathit{u211}=0,  \\
{\mathit{u12}_{\mathit{x2}}} - (\mathit{u112} + \mathit{u222})\,
\mathit{u12} - \mathit{u22}\,\mathit{u212} - \mathit{u11}\,
\mathit{u122}=0,  \\
{\mathit{u22}_{\mathit{x1}}} - 2\,\mathit{u12}\,\mathit{u112} - 2
\,\mathit{u22}\,\mathit{u212}=0, \\
{\mathit{u22}_{\mathit{x2}}}
 - 2\,\mathit{u12}\,\mathit{u122} - 2\,\mathit{u22}\,\mathit{u222
}=0,  \\
{\mathit{u111}_{\mathit{x2}}} - {\mathit{u112}_{\mathit{x1}}} - 
\mathit{u212}\,\mathit{u112} + \mathit{u122}\,\mathit{u211} \\
\qquad = \mathit{\_C1}\,\mathit{u12} + \sqrt{ - \mathit{u12}^{2} + 
\mathit{u22}\,\mathit{u11}}\,\mathit{\_C2},  \\
{\mathit{u112}_{\mathit{x2}}} - {\mathit{u122}_{\mathit{x1}}} - (
\mathit{u111} - \mathit{u212})\,\mathit{u122} - \mathit{u112}\,
\mathit{u222} + \mathit{u112}^{2}=\mathit{\_C1}\,\mathit{u22}, 
 \\
{\mathit{u211}_{\mathit{x2}}} - {\mathit{u212}_{\mathit{x1}}} - 
\mathit{u212}^{2} + \mathit{u212}\,\mathit{u111} - (\mathit{u112}
 - \mathit{u222})\,\mathit{u211}= - \mathit{u11}\,\mathit{\_C1}, 
 \\
{\mathit{u212}_{\mathit{x2}}} - {\mathit{u222}_{\mathit{x1}}} - 
\mathit{u122}\,\mathit{u211} + \mathit{u212}\,\mathit{u112} \\ 
\qquad = - \mathit{\_C1}\,\mathit{u12} + \sqrt{ - \mathit{u12}^{2} + 
\mathit{u22}\,\mathit{u11}}\,\mathit{\_C2}] }
}
\end{maplelatex}

\end{maplegroup}
\begin{maplegroup}
They show that all equivariant sections $c$ can be parametrised by two
constants ($\_C1$ and $\_C2$). The second constant $\_C2$ is special
to the $2$-dimensional case and we obtain $\_C2=0$ using the Jacobi
conditions in \cite[Thm. 7.4.8]{Pommaret78}. The first six
integrability conditions express the Christoffel symbols in terms of
the metric and its first order derivatives (cf. eq.
(\ref{eq_Christoffel_Symbol})):

\end{maplegroup}
\begin{maplegroup}
\begin{mapleinput}
\mapleinline{active}{1d}{nrsolve(Ves[1..6],uvar_Gamma)[1];}{%
}
\end{mapleinput}

\mapleresult
\begin{maplelatex}
\mapleinline{inert}{2d}{[u111 =
-1/2*(-u11[x1]*u22+2*u12*u12[x1]-u12*u11[x2])/(-u12^2+u22*u11), u112 =
-1/2*(-u11[x2]*u22+u12*u22[x1])/(-u12^2+u22*u11), u211 =
1/2*(2*u11*u12[x1]-u12*u11[x1]-u11*u11[x2])/(-u12^2+u22*u11), u212 =
1/2*(u11*u22[x1]-u12*u11[x2])/(-u12^2+u22*u11), u122 =
-1/2*(u22[x2]*u12-2*u22*u12[x2]+u22*u22[x1])/(-u12^2+u22*u11), u222 =
1/2*(u11*u22[x2]-2*u12*u12[x2]+u12*u22[x1])/(-u12^2+u22*u11)];}{%
\maplemultiline{
[\mathit{u111}= - {\displaystyle \frac {1}{2}} \,{\displaystyle 
\frac { - {\mathit{u11}_{\mathit{x1}}}\,\mathit{u22} + 2\,
\mathit{u12}\,{\mathit{u12}_{\mathit{x1}}} - \mathit{u12}\,{
\mathit{u11}_{\mathit{x2}}}}{ - \mathit{u12}^{2} + \mathit{u22}\,
\mathit{u11}}} , \\
 \mathit{u112}= - {\displaystyle \frac {1}{2}} 
\,{\displaystyle \frac { - {\mathit{u11}_{\mathit{x2}}}\,\mathit{
u22} + \mathit{u12}\,{\mathit{u22}_{\mathit{x1}}}}{ - \mathit{u12
}^{2} + \mathit{u22}\,\mathit{u11}}} ,  \\
\mathit{u211}={\displaystyle \frac {1}{2}} \,{\displaystyle 
\frac {2\,\mathit{u11}\,{\mathit{u12}_{\mathit{x1}}} - \mathit{
u12}\,{\mathit{u11}_{\mathit{x1}}} - \mathit{u11}\,{\mathit{u11}
_{\mathit{x2}}}}{ - \mathit{u12}^{2} + \mathit{u22}\,\mathit{u11}
}} , \\ 
\mathit{u212}={\displaystyle \frac {1}{2}} \,
{\displaystyle \frac {\mathit{u11}\,{\mathit{u22}_{\mathit{x1}}}
 - \mathit{u12}\,{\mathit{u11}_{\mathit{x2}}}}{ - \mathit{u12}^{2
} + \mathit{u22}\,\mathit{u11}}} ,  \\
\mathit{u122}= - {\displaystyle \frac {1}{2}} \,{\displaystyle 
\frac {{\mathit{u22}_{\mathit{x2}}}\,\mathit{u12} - 2\,\mathit{
u22}\,{\mathit{u12}_{\mathit{x2}}} + \mathit{u22}\,{\mathit{u22}
_{\mathit{x1}}}}{ - \mathit{u12}^{2} + \mathit{u22}\,\mathit{u11}
}} ,  \\
\mathit{u222}={\displaystyle \frac {1}{2}} \,{\displaystyle 
\frac {\mathit{u11}\,{\mathit{u22}_{\mathit{x2}}} - 2\,\mathit{
u12}\,{\mathit{u12}_{\mathit{x2}}} + \mathit{u12}\,{\mathit{u22}
_{\mathit{x1}}}}{ - \mathit{u12}^{2} + \mathit{u22}\,\mathit{u11}
}} ] }
}
\end{maplelatex}

\end{maplegroup}
\begin{maplegroup}
Starting with a metric, they are always fulfilled, but an arbitrary
section of $\mathcal{F}$ allows to choose metric and Christoffel
symbols independently. The last four integrability conditions express components
of the Riemann curvature tensor as derivatives of the Christoffel symbols. 
The equations are equivalent ($\_C2=0$) to the condition of a metric with constant
scalar curvature: 
\[   
R^k_{lij} = \partial_i \Gamma^k_{lj} -
\partial_j \Gamma^k_{li} + \Gamma^r_{lj} \, \Gamma^k_{ri} -
\Gamma^r_{li} \, \Gamma^k_{rj} = \_C1
      (\delta^k_j g_{li} - \delta^k_i g_{lj}).  
\]   
The calculations
with {\tt jets} complement the the theory with explicit coordinate
changes of the natural bundles $\mathcal{F}$ and $\mathcal{F}_1$,
which are equivalent to the $\Pi_q$-action on $\mathcal{F}$. Vessiot's
structure equations can now be used to check the integrability of jet
groupoids.

In \cite{Vessiot03}, the structure equations are solved for the
constants, which is an alternative choice of coordinates for $\mcF_1$. Usually,
this leads to larger expressions and hides the geometrical interpretation. If
the typical fibre of $\mathcal{F}$ is not homogeneous, the freedom for equivariant
sections may extend from constants to smooth invariants.

Vessiot's structure equations can also be applied to test whether two geometric
objects are formally equivalent, which is connected to the integrability of the
corresponding groupoids an equivariant sections.

\end{maplegroup}
\begin{maplegroup}
\begin{mapleinput}
\end{mapleinput}

\end{maplegroup}
\begin{maplegroup}
\begin{mapleinput}
\end{mapleinput}

\end{maplegroup}

\nocite{Sternberg64}

\bibliography{literature}

\def\cprime{$'$}
\providecommand{\bysame}{\leavevmode\hbox to3em{\hrulefill}\thinspace}
\providecommand{\MR}{\relax\ifhmode\unskip\space\fi MR }
\providecommand{\MRhref}[2]{%
  \href{http://www.ams.org/mathscinet-getitem?mr=#1}{#2}
}
\providecommand{\href}[2]{#2}
\begin{thebibliography}{10}

\bibitem{Barakat01}
Mohamed Barakat, \emph{Jets. {A} {MAPLE}-package for formal differential
  geometry}, Compu\-ter algebra in scientific computing (Konstanz, 2001),
  Springer, Berlin, 2001, pp.~1--12. \MR{MR1942047 (2003j:53001)}

\bibitem{Ehresmann55}
Charles Ehresmann, \emph{Les prolongements d'un espace fibr\'e
  diff\'erentiable}, C. R. Acad. Sci. Paris \textbf{240} (1955), 1755--1757.
  \MR{MR0071083 (17,80c)}

\bibitem{Goldschmidt67DG}
Hubert Goldschmidt, \emph{Integrability criteria for systems of nonlinear
  partial differential equations}, J. Differential Geometry \textbf{1} (1967),
  269--307. \MR{MR0226156 (37 \#1746)}

\bibitem{KMS93}
Ivan Kol{\'a}{\v{r}}, Peter~W. Michor, and Jan Slov{\'a}k, \emph{Natural
  operations in differential geometry}, Springer-Verlag, Berlin, 1993.
  \MR{MR1202431 (94a:58004)}

\bibitem{KLV86}
I.~S. Krasil{\cprime}shchik, V.~V. Lychagin, and A.~M. Vinogradov,
  \emph{Geometry of jet spaces and nonlinear partial differential equations},
  Advanced Studies in Contemporary Mathematics, vol.~1, Gordon and Breach
  Science Publishers, New York, 1986, Translated from the Russian by A. B.
  Sosinski\u\i. \MR{MR861121 (88m:58211)}

\bibitem{Lie1891}
S.~Lie, \emph{{D}ie {G}rundlagen f\"ur die {T}heorie der unendlichen
  kontinuierlichen {T}ransformationsgruppen {I}/{II}}, Leipz. Ber. \textbf{III}
  (1891), 316--393.

\bibitem{Mackenzie05}
Kirill C.~H. Mackenzie, \emph{General theory of {L}ie groupoids and {L}ie
  algebroids}, London Mathematical Society Lecture Note Series, vol. 213,
  Cambridge University Press, Cambridge, 2005. \MR{MR2157566 (2006k:58035)}

\bibitem{MoerMrcun03}
I.~Moerdijk and J.~Mr{\v{c}}un, \emph{Introduction to foliations and {L}ie
  groupoids}, Cambridge Studies in Advanced Mathematics, vol.~91, Cambridge
  University Press, Cambridge, 2003. \MR{MR2012261 (2005c:58039)}

\bibitem{Pommaret78}
J.-F. Pommaret, \emph{Systems of partial differential equations and {L}ie
  pseudogroups}, Mathematics and its Applications, vol.~14, Gordon \& Breach
  Science Publishers, New York, 1978, With a preface by Andr\'e Lichnerowicz.
  \MR{MR517402 (81f:58046)}

\bibitem{Spencer69}
D.~C. Spencer, \emph{Overdetermined systems of linear partial differential
  equations}, Bull. Amer. Math. Soc. \textbf{75} (1969), 179--239.
  \MR{MR0242200 (39 \#3533)}

\bibitem{Sternberg64}
Shlomo Sternberg, \emph{Lectures on differential geometry}, Prentice-Hall Inc.,
  Englewood Cliffs, N.J., 1964. \MR{MR0193578 (33 \#1797)}

\bibitem{Terng78}
Chuu~Lian Terng, \emph{Natural vector bundles and natural differential
  operators}, Amer. J. Math. \textbf{100} (1978), no.~4, 775--828. \MR{MR509074
  (81c:58009)}

\bibitem{Vessiot03}
Ernest Vessiot, \emph{Sur la th\'eorie des groupes continus}, Ann. Sci. \'Ecole
  Norm. Sup. (3) \textbf{20} (1903), 411--451. \MR{MR1509031}

\end{thebibliography}
\bibliographystyle{amsplain}

\end{document}